\documentclass[twoside,leqno]{article}
\usepackage[letterpaper]{geometry}
\usepackage{siamproceedings}
\usepackage{amsfonts,amssymb,amscd,mathtools}
\usepackage{tikz} 
\usetikzlibrary{decorations.pathreplacing}
\newsiamremark{remark}{Remark}
\newsiamthm{conjecture}{Conjecture}
\newsiamremark{example}{Example}
\def\R{\mathbb{R}}
\def\C{\mathbb{C}}
\def\Z{\mathbb{Z}}
\def\K{\mathbb{K}}
\def\Q{\mathbb{Q}}
\def\O{\mathcal{O}}
\def\CP{\mathbb{CP}}
\def\PP{\mathbb{P}}
\def\F{\mathcal{F}}
\def\W{\mathcal{W}}
\def\Coh{\mathrm{Coh}}
\def\bb{\mathfrak{b}}
\def\hol{\mathrm{hol}}

\begin{document}
\title{Lagrangian Floer theory, from geometry to algebra and back again}
\author{Denis Auroux\thanks{Department of Mathematics, Harvard University,
Cambridge MA 02138, USA (\email{auroux@math.harvard.edu}).}}
\date{}

\maketitle


\begin{abstract} 
We survey various aspects of Floer theory and its place in modern symplectic geometry,
from its introduction to address classical conjectures of Arnold about
Hamiltonian diffeomorphisms and Lagrangian submanifolds, to the rich algebraic 
structures captured by the Fukaya category, and finally to the idea,
motivated by mirror symmetry, of a ``geometry of Floer theory'' centered
around family Floer cohomology and local-to-global principles for
Fukaya categories.
\end{abstract}

\section{Arnold's conjectures and Floer (co)homology.}

\ Since their introduction by Andreas Floer almost forty years ago to study
fundamental questions in symplectic geometry and low-dimensional topology,
Floer homology theories have been a major driver of progress in those areas
of mathematics.

Floer's work in symplectic geometry was motivated by key conjectures of
Arnold about {\em fixed points of Hamiltonian diffeomorphisms} and
{\em intersections of Lagrangian submanifolds} \cite{Arnold}. 
We briefly recall some basic
notions in symplectic geometry (see \cite{Cannas} for a detailed treatment).
First of all, a {\bf symplectic form} is a closed non-degenerate 2-form
$\omega$ on a smooth manifold $M$. For instance:
\begin{itemize}
\item $\R^{2n}=\C^n$ carries the standard symplectic
form $\omega_0=\sum dx_i\wedge dy_i$;
\item any cotangent bundle $M=T^*N$ has an exact symplectic form 
$\omega=d\lambda$, given in coordinates $(q_i)$ on $N$ and dual
coordinates $(p_i)$ on the fibers by $\omega=\sum dp_i\wedge
dq_i$; here $\lambda=\sum p_i\,dq_i$ is the canonical Liouville form;
\item $\CP^n$, or by restriction, any complex projective variety, with the Fubini-Study K\"ahler
form.
\end{itemize}
A smooth function $H:M\to \R$ on a symplectic manifold $(M,\omega)$
determines a {\bf Hamiltonian vector field} $X_H$, characterized by the
property that $\omega(\cdot,X_H)=dH$.  The flows generated by time-dependent 
Hamiltonian vector fields are called {\bf Hamiltonian diffeomorphisms};
they form a subgroup $\mathrm{Ham}(M,\omega)$ of the 
symplectomorphism group $\mathrm{Symp}(M,\omega)=\{\varphi\in
\mathrm{Diff}(M)\,|\,\varphi^*\omega=\omega\}$. For instance, the classical
mechanical system consisting of a point mass moving in a potential $V(q)$
on the manifold $N$ is described by the dynamics generated by the Hamiltonian $H(q,p)=\frac12 |p|^2+V(q)$
on the phase space $T^*N$. 
\smallskip

Arnold's conjecture states that
the number of fixed points of a Hamiltonian diffeomorphism of a compact
symplectic manifold $(M,\omega)$ is at least the minimal number of critical 
points of a smooth function on $M$. Moreover:

\begin{conjecture}[\bf Arnold's conjecture for non-degenerate
Hamiltonians]\label{conj:arnoldconj}\!
Given a compact symplectic manifold $(M,\omega)$ and $\varphi\in
\mathrm{Ham}(M,\omega)$ with non-degenerate fixed points,\vspace*{-4pt}
\begin{equation}\label{eq:arnoldbound}
\#\mathrm{Fix}(\varphi)\geq \sum_{i=0}^{\dim M} \dim H^i(M,\Q).
\end{equation}
\end{conjecture}

Arnold's conjecture in this formulation was first proved by 
Floer for monotone symplectic manifolds \cite{Floer,Floer2}, then
in the semi-positive setting by
Hofer-Salamon \cite{HoferSalamon} and Ono \cite{Ono}, and in full generality
by  Fukaya-Ono \cite{FOno} and Liu-Tian \cite{LiuTian} (see also Pardon \cite{Pardon} or 
Filippenko-Wehrheim \cite{FW} for more modern treatments).

The inequality \eqref{eq:arnoldbound} is easily seen to hold for Hamiltonian diffeomorphisms
generated by a time-independent Morse function, since $X_H$ vanishes at every critical point of $H$. 
The general case however is considerably more difficult, and the proof 
crucially relies on the construction of Hamiltonian Floer
(co)homology,\footnote{As with Morse
(co)homology, the only differences between Floer homology and cohomology 
are in the grading conventions and in deciding which end of a trajectory is
the output. We will consistently prefer cohomology, due to its ring
structure.} which we now sketch. (See \cite{AudinDamian} and
\cite[Chapter 12]{McS} for detailed texts about the construction.) 

\subsection{\!\!Hamiltonian Floer (co)homology and Arnold's conjecture.}

Given a (time-dependent, 1-periodic, non-degenerate) Hamiltonian $H$ and an $\omega$-compatible
almost-complex structure $J$ on $M$ (i.e., $J\in \mathrm{End}(TM)$
satisfying $J^2=-1$ and such that $\omega(\cdot,J\cdot)$ is a Riemannian
metric), the Floer complex $CF^*(M,H;J)$ is the free module (over a suitable
coefficient ring or field) generated by the set $\mathcal{X}(H)$ of 1-periodic orbits of
$X_H$, equipped
with a differential which counts solutions to a perturbed Cauchy-Riemann
equation known as {\bf Floer's equation}. Namely, we consider maps
$u:\R\times S^1\to M$, $(s,t)\mapsto u(s,t)$, which satisfy \vspace*{-1mm}
\begin{equation}
\label{eq:floereq}
\frac{\partial u}{\partial s}+J\Bigl(\frac{\partial u}{\partial t}-X_{H_t}\Bigr)=0,%
\vspace*{-1mm}
\end{equation}
with finite energy $E(u)=\iint \|\partial u/\partial s\|^2\,ds\,dt$,
and which converge as $s\to \pm \infty$ to given periodic orbits $x_{\pm}\in
\mathcal{X}(H)$.
Floer's equation is (formally) the gradient flow of the
{\bf action functional} on the free loop space of $M$, defined by%
\vspace*{-1mm}
\begin{equation}\label{eq:action}
\mathcal{A}_H(x)=\int \bigl(
H_t(x(t))-\lambda(\dot x(t))\bigr)\,dt
\vspace*{-1mm}
\end{equation}
when $\omega=d\lambda$ is exact (and in this case
$E(u)=\mathcal{A}_H(x_-)-\mathcal{A}_H(x_+)$); in the non-exact case $\mathcal{A}_H$ is multivalued, but the intuition
remains the same: the 1-periodic orbits of $X_H$ are the critical points of 
$\mathcal{A}_H$, and the Floer (co)homology of $H$ is (formally) the
Morse(-Novikov) (co)homology of $\mathcal{A}_H$.

For generic $H$ and $J$, the moduli space $\mathcal{M}(x_+,x_-;[u],J)$
of Floer trajectories (i.e., solutions to \eqref{eq:floereq})
connecting given generators $x_\pm\in \mathcal{X}(H)$ and in a given homotopy class
$[u]$,
up to $\R$-translation in the $s$ direction, is a finite dimensional
manifold (of dimension determined by the Conley-Zehnder index). 

As an application of Gromov compactness, $\mathcal{M}(x_+,x_-;[u],J)$ can be
compactified by adding (1) {\bf broken trajectories}, i.e.\ 
tuples of Floer trajectories from $x_+$ to some other generator, and so on, to
$x_-$, and (2) nodal configurations consisting of Floer trajectories (of lower energy) 
together with one or more $J$-holomorphic {\bf sphere bubbles} attached to each other at nodes.

When all the strata of the compactified moduli space
$\overline{\mathcal{M}}(x_+,x_-;[u],J)$ are regular (i.e., smooth and of the expected dimension),
sphere bubbling only occurs in real codimension 2, and the (codimension 1) boundary
corresponds exactly to two-step broken trajectories:
\vspace*{-1.5mm}
\begin{equation}\label{eq:brokentraj}
\partial \overline{\mathcal{M}}(x_+,x_-;[u],J)=\bigsqcup_{\substack{y\in
\mathcal{X}(H)\\ [u]=[u_-]\#[u_+]}} \overline{\mathcal{M}}(y,x_-;[u_-],J)\times
\overline{\mathcal{M}}(x_+,y;[u_+],J).
\end{equation}
The {\bf Floer differential} $\partial:CF^*(M,H;J)\to CF^*(M,H;J)$ is then defined by
counting trajectories in zero-dimensional moduli spaces, with signs and with suitable
weights $w([u])$:
\vspace*{-1mm}
\begin{equation}\label{eq:floerdiff}
\partial x_+=\sum_{x_-,[u]} \#(\mathcal{M}(x_+,x_-;[u],J))\,w([u])\,x_-.
\end{equation}\vskip-1.5mm

\noindent (In the non-exact setting, one typically works over the {\bf Novikov field} of formal
power series with real exponents,
$\mathbb{K}=\Lambda_k:=\bigl\{\sum a_i T^{\lambda_i},\ a_i\in k,\ \lambda_i\in\R,\
\lambda_i\to +\infty\bigr\}$, and solutions are weighted by energy: $w([u])=T^{E([u])}$.)
As in classical Morse theory, by considering \eqref{eq:brokentraj} for 1-dimensional moduli
spaces one finds that $\partial^2=0$, which allows one to define the Floer
cohomology $HF^*(M,H;J)=\mathrm{Ker}(\partial)/\mathrm{Im}(\partial)$.

When $M$ is compact, the Floer cohomology is independent of the choice of
compatible almost-complex structure $J$ and Hamiltonian $H$, and it
is isomorphic to the (quantum) cohomology of $M$ via the PSS (Piunikhin-Salamon-Schwarz)
isomorphism. (This is in fact a ring isomorphism, where the product structure on
Floer cohomology is defined by counting solutions of a suitable analogue of
\eqref{eq:floereq} for maps from a pair of pants to $M$.) Thus, when
Hamiltonian Floer cohomology is well-defined, the Arnold conjecture follows
from the fact that the rank of the Floer complex is bounded from below by the rank of its
cohomology. A key technical point that restricted the level of generality of
the early proofs, however, is achieving regularity for moduli spaces of
configurations that include sphere bubbles so that \eqref{eq:brokentraj}
holds (see \cite{McS} for the argument
in the semi-positive case); or, failing that, constructing {\bf virtual fundamental classes} for the moduli spaces
\cite{FOno,LiuTian,Pardon,FO3kuranishi}. 

\begin{remark}
Outside of the monotone setting, the construction of Floer cohomology usually requires working over a
field of characteristic zero (due to the moduli spaces being orbifolds 
rather than manifolds in the presence of multiply covered bubbles), 
thus yielding the bound \eqref{eq:arnoldbound} by the cohomology of $M$ with
rational coefficients. Two recent advances make it possible to improve on
this bound when the cohomology of $M$ has torsion. First, Abouzaid and Blumberg 
have shown that Floer moduli spaces admit fundamental classes
not only in ordinary cohomology but also in Morava K-theory \cite{AbBlumberg}, which
leads to a bound by the rank of the cohomology of $M$ with
coefficients in a field of characteristic $p$. Second, Bai-Xu \cite{BaiXu}
and Rezchikov \cite{Rezchikov} have constructed a version of Floer theory
over the Novikov ring with $\Z$ coefficients (by dealing with
multiply covered configurations in a specific manner), which gives a lower
bound by the total Betti numbers over $\Z$ (including torsions of all
characteristics). 
\end{remark}

Floer cohomology can also be defined on noncompact symplectic manifolds, 
subject to suitable conditions on the manifold and on the Floer data $(H,J)$, 
so as to control the behavior of solutions to \eqref{eq:floereq} near
infinity.  A common setup is that of {\bf Liouville manifolds},
i.e.\ exact symplectic manifolds whose ends are modelled on the
positive symplectization of a contact manifold $(N,\alpha)$, namely $N\times
[1,+\infty)$ with the symplectic form $\omega=d(r\alpha)$.
(For instance, cotangent bundles are Liouville, and so are affine complex algebraic
varieties.) One then considers Hamiltonians which grow
linearly with $r$ at infinity (whose 1-periodic orbits near infinity 
correspond to closed orbits of the Reeb vector field of the contact
form $\alpha$ on $N$); the Floer cohomology of such $H$ depends substantially on
its slope at infinity, but the direct limit of the Floer
cohomology groups for Hamiltonians of increasingly large slope, called {\bf
symplectic cohomology} and denoted $SH^*(M)$, is independent of auxiliary choices.

\subsection{\!\!Lagrangian Floer (co)homology.}

Lagrangian Floer theory was developed by Floer \cite{Floer} to approach
another conjecture of Arnold about intersections of Lagrangian submanifolds. 
Recall that an $n$-dimensional submanifold of a $2n$-dimensional 
symplectic manifold $(M,\omega)$ is said to be {\bf Lagrangian} if
$\omega_{|L}=0$. For example, the zero section is a Lagrangian
submanifold of the cotangent bundle $T^*N$ with its standard symplectic
form, and in fact, by the Weinstein neighborhood theorem, this is the universal
example: a Lagrangian submanifold $L\subset M$ always admits a tubular neighborhood
which is symplectomorphic to a neighborhood of the zero section in $T^*L$.

The {\bf Arnold-Givental conjecture} and other Lagrangian versions of
Arnold's conjecture state that, under suitable assumptions, the image $\varphi(L)$ of a
Lagrangian submanifold $L$ under a Hamiltonian diffeomorphism $\varphi\in
\mathrm{Ham}(M,\omega)$, when it intersects $L$ transversely, must do so in
a number of points at least equal to the total dimension of the cohomology
of $L$. The statement doesn't hold unconditionally: for example the unit
circle $S^1\subset (\R^2,dx\wedge dy)$ is disjoint from its image under the
translation generated by a sufficiently large linear Hamiltonian. On the
other hand, it holds for the diagonal Lagrangian $\Delta\subset (M\times
M,-\omega\oplus \omega)$. Since the image of $\Delta$ under a Hamiltonian diffeomorphism of
the form $\mathrm{id}\times \varphi$ is the graph of $\varphi$, the
Arnold conjecture for Hamiltonians follows as a special case.
Floer proved:

\begin{theorem}[Floer \cite{Floer}]\label{thm:floerlag}
Let $L$ be a compact Lagrangian submanifold of a compact symplectic
manifold $(M,\omega)$, such that $\pi_2(M,L)=0$. Then for any $\varphi\in
\mathrm{Ham}(M,\omega)$ such that $\varphi(L)$ meets $L$ transversely,
\vskip-2mm
$$\#(L\cap \varphi(L))\geq \sum_{i=0}^{\dim L} \dim H^i(L,\Z/2).$$
\end{theorem}
The assumption can be relaxed to only require $[\omega]\cdot
\pi_2(M,L)=0$, i.e. the integral of
$\omega$ over any disc bounded by $L$ is zero; and $M$ need not be compact
(but should be convex at infinity, so Floer theory is well-defined). Thus, Theorem
\ref{thm:floerlag} applies e.g.\ to a compact {\bf exact Lagrangian}
submanifold of a Liouville manifold $(M,\omega=d\lambda)$, i.e.\ a Lagrangian
submanifold $L$ such that the 1-form $\lambda_{|L}$ is exact
(e.g., the zero section in a cotangent bundle).
\smallskip

To a pair of Lagrangian submanifolds $L_0,L_1\subset
(M,\omega)$, a (time-dependent) Hamiltonian $H$ such that $\varphi_H^1(L_0)$
intersects $L_1$ transversely, and an $\omega$-compatible
almost complex structure $J$ on $M$, one associates the {\bf Lagrangian Floer complex}
$CF^*(L_0,L_1;H,J)$, which is a free module (over a suitable coefficient
field $\mathbb{K}$:
$\Z/2$ in the setting of Theorem \ref{thm:floerlag}, the Novikov field in
non-exact settings,
etc.) generated by the set $\mathcal{X}=\mathcal{X}(L_0,L_1;H)$ of
trajectories $x(t)$ of the Hamiltonian vector field $X_H$ such that $x(0)\in L_0$ 
and $x(1)\in L_1$. (Equivalently, since each trajectory is determined by its
end point $x(1)$, one could also define 
$\mathcal{X}=\varphi_H^1(L_0)\cap L_1$.) This is equipped with the Floer
differential, which counts solutions to Floer's equation \eqref{eq:floereq},
where now the domain of the map $u$ is the strip $\R\times [0,1]$, and we
impose boundary conditions $u(s,0)\in L_0$ and $u(s,1)\in L_1$ for all $s\in
\R$. (As in the Hamiltonian case, solutions to \eqref{eq:floereq} can also
be thought of as gradient flow lines of a suitably defined action functional
on the space of paths $x:[0,1]\to M$ which start on $L_0$ and end on $L_1$.)

For generic $H$ and $J$, the moduli space $\mathcal{M}(x_+,x_-;[u],J)$
of Floer trajectories connecting given generators $x_\pm \in \mathcal{X}$
and in a given homotopy class $[u]$, up to $\R$-translation, is a finite
dimensional manifold (of dimension determined by the Maslov index).
As in the Hamiltonian case, the Floer differential
$\partial:CF^*(L_0,L_1;H,J)\to CF^*(L_0,L_1;H,J)$ is defined by
counting trajectories in zero-dimensional moduli spaces (with suitable
signs and weights $w([u])$ when the coefficient field is not $\Z/2$ as in
the setting of Theorem \ref{thm:floerlag}):
\begin{equation}\label{eq:floerdiff2}
\partial x_+=\sum_{x_-,[u]} \#(\mathcal{M}(x_+,x_-;[u],J))\,w([u])\,x_-.
\end{equation}\vskip-1mm
\noindent
(When $\mathrm{char}(\mathbb{K})\neq 2$, one should equip $L_0$ and $L_1$ with spin
structures in order to determine orientations of the moduli spaces so as
to be able to count with signs.)

The compactification of $\mathcal{M}(x_+,x_-;[u],J)$ involves
adding broken trajectories (which arise as limits of configurations
where energy escapes towards $s\to \pm \infty$) as well as nodal configurations
containing $J$-holomorphic sphere or disc bubbles (which arise as limits when energy
concentrates near an interior or boundary point of the domain
$\R\times [0,1]$). 
In general, {\bf disc bubbling} is an
{\bf obstruction} to the definition of Lagrangian Floer cohomology, as it
is expected to occur in real codimension 1, thereby causing the boundary of
the compactified moduli space $\overline{\mathcal{M}}(x_+,x_-;[u],J)$ to
{\em not} satisfy \eqref{eq:brokentraj}. This in turn means that the
Floer differential defined by \eqref{eq:floerdiff2} does not always
square to zero. Thus, excluding the occurrence of disc bubbles (or
cancelling them out algebraically) is a prerequisite to the definition of
Lagrangian Floer cohomology. (As in the Hamiltonian case, 
there are also issues of regularity of the moduli spaces, or how to 
count solutions in the absence of regularity, to be dealt with. In simple
cases these can be dealt with by elementary methods, but the general case
requires abstract perturbation techniques such as Kuranishi structures or
polyfolds; we refer the reader to \cite{FO3kuranishi} for
one possible treatment.)

In the setting of Theorem \ref{thm:floerlag}, the assumption that $[\omega]\cdot \pi_2(M,L_i)=0$ ensures
that $L_i$ cannot bound any holomorphic discs (since  $J$-holomorphic curves have positive
symplectic area), thereby excluding disc
bubbling (as well as sphere bubbling, since it follows that $[\omega]\cdot
\pi_2(M)=0$) and making Floer cohomology well-defined.

For compact Lagrangians, the Floer complexes for different choices of Floer data
(i.e., the Hamiltonian $H$ and the almost-complex structure $J$) are related
to each other by {\em continuation maps} which count solutions to a version of
Floer's equation \eqref{eq:floereq} where $J$ and $H$ depend on $s$. 
In the absence of disc bubbling, the continuation maps are
quasi-isomorphisms, so that the Floer cohomology $HF^*(L_0,L_1)$ is
independent of these auxiliary choices. Theorem~\ref{thm:floerlag} then
follows from an explicit comparison between the Floer complex
$CF^*(L,L;H,J)$ for a Hamiltonian $H$ whose restriction to $L$ is a small
multiple of a given Morse function $f:L\to\R$, and the Morse complex of $f$,
from which one deduces that $HF^*(L,L)\simeq H^*(L;\mathbb{K})$ under the assumption
that $L$ is compact and does not bound any $J$-holomorphic discs.

In the presence of $J$-holomorphic discs, there are situations where disc
bubbling can be kept under control and Floer cohomology can still be defined
by elementary methods. One of these is the case of monotone Lagrangian
submanifolds in monotone symplectic manifolds (i.e., when the
Maslov class and the symplectic area of disks are positively proportional to
each other), first studied by Oh \cite{Oh}. 
However, the definition of Lagrangian Floer cohomology in full generality requires
both addressing the regularity issues that may arise with moduli spaces, and
coming up with a way to study and, when possible, cancel out the
obstructions that come from disc bubbling. This large-scale undertaking is
the heart of Fukaya, Oh, Ohta and Ono's monograph \cite{FO3book}.
Fukaya-Oh-Ohta-Ono deal with the latter issue by introducing the notion of 
{\em bounding cochain}, which we will discuss in the next section. This has
the following striking application to the Lagrangian Arnold conjecture:

\begin{theorem}[{Fukaya-Oh-Ohta-Ono \cite[Theorem H]{FO3book}}]\label{thm:floerlagFO3}
Let $L$ be a relatively spin\footnote{A submanifold $L\subset M$ is said to
be relatively spin if $w_2(TL)$ is in the image
of the restriction map $i^*:H^2(M,\Z/2)\to H^2(L,\Z/2)$.} compact Lagrangian 
submanifold of a compact symplectic
manifold $(M,\omega)$. Assume that the inclusion map $i_*:H_*(L,\Q)\to
H_*(M,\Q)$ is injective. Then for any $\varphi\in
\mathrm{Ham}(M,\omega)$ such that $\varphi(L)$ meets $L$ transversely,
\vskip-2mm
$$\#(L\cap \varphi(L))\geq \sum_{i=0}^{\dim L} \dim H_i(L,\Q).$$
\end{theorem}

\noindent The key ingredient that goes into the proof of this result is that
the injectivity of $i_*$ ensures the existence of a 
bounding cochain $b\in CF^*(L,L)$ which can be used to algebraically cancel 
the effects of disc bubbling.\footnote{More precisely:
there exists a {\em bulk deformation} of Floer theory for which $L$
admits a {\em weak bounding cochain}; see \S \ref{ss:floerop} and
\S \ref{ss:ocdeforms}.}
\medskip

Floer cohomology can also be defined for noncompact Lagrangian submanifolds
in noncompact symplectic manifolds, subject to suitable conditions on the
geometry at infinity. For example, one can consider exact Lagrangian submanifolds
of a Liouville manifold which are cylindrical at infinity, i.e.\ whose ends
are modelled on cylinders $\Lambda\times [R,+\infty)$ in the
symplectization $(N\times [1,+\infty),\omega=d(r\alpha))$, where $\Lambda$
is a {\em Legendrian} submanifold of the contact manifold $(N,\alpha)$,
i.e.\ an $(n-1)$-dimensional submanifold such that $\alpha_{|\Lambda}=0$.
Using Hamiltonians which grow linearly with $r$ at infinity, and taking the
direct limit of the Floer cohomology groups associated to Hamiltonians of
increasingly large slopes, one arrives at the {\bf wrapped Floer
cohomology} $HW^*(L_0,L_1)$ \cite{AS}; geometrically, the generators of the wrapped
Floer complex consist of intersections of $L_0$ with $L_1$ in the interior
of the manifold as well as trajectories of $X_H$ from $L_0$ to $L_1$ near
infinity, which correspond to {\em Reeb chords} between the Legendrians 
$\Lambda_0$ and $\Lambda_1$ on which the ends of $L_0$ and $L_1$ are
modelled. One can also consider other classes of Hamiltonians, giving rise
to so-called ``partially wrapped'' Floer cohomologies; there are 
natural continuation maps from ``less wrapped'' to ``more wrapped'' flavors
of Floer cohomology, but they need not be isomorphisms.
\medskip

We conclude this section with two remarks that highlight
the central role of Lagrangian Floer cohomology among Floer theories.

\begin{remark}
Hamiltonian Floer cohomology is a special case of Lagrangian
Floer cohomology. Given a symplectic manifold $(M,\omega)$, the diagonal
$\Delta_M$ is a Lagrangian submanifold in the product $M\times M$ equipped
with the symplectic form $-\pi_1^*\omega+\pi_2^*\omega$; and given $\varphi
\in \mathrm{Ham}(M,\omega)$, the Hamiltonian diffeomorphism
$\mathrm{id}\times\varphi$ of $M\times M$ maps $\Delta_M$ to the graph of
$\varphi$. There is a natural isomorphism between the Floer complexes 
$CF^*(M,H)$ and $CF^*(\Delta_M,\Delta_M)$ (for suitable choices of
Hamiltonian and almost-complex structure on $M\times M$) which
intertwines the Floer differentials, so that the Floer cohomologies are
isomorphic as well. In fact, Arnold's bound \eqref{eq:arnoldbound} for fixed points of 
$\varphi$ is equivalent to the corresponding bound for Lagrangian intersections
of $\Delta_M$ and its image under $\mathrm{id}\times\varphi$, and so
Conjecture \ref{conj:arnoldconj} also follows from Theorem
\ref{thm:floerlagFO3}.
\end{remark}

\begin{remark}\label{rmk:lowdimtop}
In a different direction, low-dimensional topology has greatly benefitted
from the introduction and study of various Floer-type invariants of 3-manifolds
and knots and links in them. 
Among these, Ozsv\'ath and Szab\'o's {\bf Heegaard-Floer homology} \cite{OS} is
most directly related to Lagrangian Floer theory: starting from a Heegaard
splitting of a 3-manifold $Y$ into two genus $g$ handlebodies glued to each 
other along their boundary surface $\Sigma$, the Heegaard-Floer homology of
$Y$ is essentially the Lagrangian Floer cohomology of a pair of Lagrangian
tori (associated to the two handlebodies) in the symmetric product
$\mathrm{Sym}^g(\Sigma)$ (or rather, in the complement of the divisor
$\{z\}\times \mathrm{Sym}^{g-1}(\Sigma)$ for the invariant 
$\widehat{HF}(Y)$, or using coefficient weights that keep track of
intersection numbers with that divisor for the invariants $HF^\pm(Y)$).
Other Floer invariants of 3-manifolds (instanton Floer homology, monopole
Floer homology) are gauge-theoretic in nature (i.e., they consider PDEs
involving connections on certain bundles over 3-manifolds and their products
with $\R$); these are related to Lagrangian Floer theory via the {\bf Atiyah-Floer
conjecture} \cite{Atiyah}, which, given a Heegaard splitting of $Y$ as above, 
relates the instanton Floer 
homology of $Y$ to the Lagrangian Floer homology of a pair of 
Lagrangian submanifolds in a moduli space
of flat connections over the surface $\Sigma$. See e.g.\ the work of
Salamon-Wehrheim \cite{SalamonWehrheim} and Daemi-Fukaya-Lipyanskiy
\cite{DFL} for progress towards the conjecture.
In fact, even invariants that aren't obviously Floer-theoretic in nature
have been shown to admit interpretations in terms of Lagrangian Floer
theory: for example, Khovanov homology \cite{AbouzaidSmithKhovanov}, or
potentially a broader range of knot homology theories \cite{Aganagic}.
\end{remark}

\section{The algebra of Floer theory: the Fukaya category}

\subsection{\!\!Operations on Floer complexes.}\label{ss:floerop}

While Floer theory on its own is already immensely useful, its full power
comes from the rich algebraic structures it carries. These are defined by
counts of solutions to Floer's equation on more general domains than the
strips and cylinders encountered in the previous section. 

Given a
Riemann surface $S$ with boundary and punctures, a 1-form on
$S$ with values in Hamiltonian vector fields on $M$, i.e.\ 
$X_H\otimes \beta$ where $\beta\in \Omega^1(S,\R)$ and $H\in
C^\infty(S\times M,\R)$, and an almost-complex structure $J$ on $M$
(possibly varying over $S$ as well), one can consider the moduli space of solutions to Floer's equation 
\begin{equation}\label{eq:floereqS}
(du-X_H\otimes \beta)^{0,1}_J=0,
\end{equation}
with Lagrangian boundary conditions along $\partial S$, of finite energy
\vspace*{-1mm}
\begin{equation}\label{eq:energy}
E(u)=\bigl\|du-X_H\otimes \beta\bigr\|^2_{L^2}=\int_S u^*\omega-u^*dH\wedge \beta,
\vspace*{-1mm}
\end{equation}
and asymptotic to
given trajectories of $X_H$ near the punctures of $S$. Near the
punctures the Hamiltonian perturbation $X_H\otimes \beta$ is required 
to be consistent with the choices made in the definition of Lagrangian or 
Hamiltonian Floer cohomology, i.e.\ of the
form $X_{H_t}\otimes dt$ in local coordinates $z=s+it$ identifying each end of $S$ 
near a boundary (resp.\ interior) puncture with a semi-infinite strip
(resp.\ cylinder), so that \eqref{eq:floereqS} reduces to \eqref{eq:floereq}.

Letting the choice of domain $S$ vary over a suitable moduli space of
Riemann surfaces with fixed topology, and considering solutions of Floer's
equation which converge to given trajectories of $X_H$ at the
punctures of $S$, we arrive at moduli spaces of
solutions to \eqref{eq:floereqS}, which can be used to 
define operations on Floer complexes via weighted counts of
solutions in moduli spaces of expected dimension zero.
(The inputs of the operation correspond to the ``positive'' strip-like or cylindrical
ends, i.e.\ those in which the local coordinate $s$ goes to $+\infty$, while the
outputs correspond to the ``negative'' ends, those at which $s\to -\infty$.)

These moduli spaces can be
compactified by allowing for broken configurations (i.e., letting the domain $S$
degenerate to a disjoint union of Riemann surfaces, 
equipped with matching pairs of strip-like or cylindrical ends at which the
solutions to Floer's equation on the various components converge to
the same trajectories of $X_H$), as well as degenerations of the
domain and/or $J$-holomorphic disc or sphere bubbling. The codimension 1 boundary strata
of the compactified moduli spaces typically consist of two-component
configurations; if one chooses the Floer perturbation data in a manner that
behaves consistently with respect to these degenerations, the boundary can
then be identified with a union of products of simpler moduli spaces. By considering the
0-dimensional boundaries of 1-dimensional moduli spaces, one arrives at algebraic relations
satisfied by the operations of the ``{open-closed Floer TQFT}''. See e.g.\ \cite[chapter
8]{SeBook} for more details.

In the setting of Hamiltonian Floer cohomology and symplectic cohomology,
the most important operation is the {\bf pair-of-pants product}, which
is defined by counting solutions to Floer's equation on a 3-punctured
Riemann sphere (i.e., a pair of pants), with two of the punctures corresponding to
inputs of the product operation and the third one to its output. The chain-level product satisfies
the Leibniz rule with respect to the Floer differential, and the induced
product on Floer cohomology (or symplectic cohomology) is associative and commutative
(in the graded sense). Moreover, the PSS construction (or an
analysis of the limit for $C^2$-small time-independent Hamiltonians) shows that (over suitable
coefficients) the Floer cohomology ring
of a compact symplectic manifold is isomorphic to its quantum cohomology.

Likewise, there is a product operation on Lagrangian Floer cohomology, which
comes from counting solutions to Floer's equation on a disc with three
boundary punctures. Counting (with appropriate weights) such
discs whose boundary arcs map to given Lagrangian submanifolds $L_0,L_1,L_2$
and whose boundary punctures converge to generators of the Floer complexes
$CF^*(L_0,L_1)$, $CF^*(L_1,L_2)$ at the inputs and $CF^*(L_0,L_2)$ at the
output (we omit the Hamiltonian perturbations from the notation), 
one obtains the chain-level product
$$\mu^2:CF^*(L_1,L_2)\otimes CF^*(L_0,L_1)\to CF^*(L_0,L_2),$$
which in the absence of disc bubbling satisfies the Leibniz rule with respect to the differentials
and defines a unital associative product on the Floer cohomology groups.

The chain-level product is only associative up to homotopy: the
Floer differential $\partial=\mu^1$ and the product $\mu^2$ are part of a
sequence of {\bf $A_\infty$-operations} (higher products) on Lagrangian Floer complexes,
\begin{equation}\label{eq:mu_d}\mu^d:CF^*(L_{d-1},L_d)\otimes \dots \otimes CF^*(L_0,L_1)\to
CF^*(L_0,L_d)\end{equation}
for $d\geq 0$, counting solutions to Floer's equation on a disc with $d+1$
boundary punctures ($d$ inputs and one output), with the successive boundary
arcs mapping to $L_0,L_1,\dots,L_d$. (When the Floer complexes
are $\Z$-graded, e.g.\ when the first Chern class $c_1(TM)$ vanishes and the
Lagrangian submanifolds $L_0,\dots,L_d$ have vanishing Maslov class, the
operation $\mu^d$ has degree $2-d$, reflecting the fact that the moduli space of
conformal structures on a $d+1$-pointed disc, which compactifies to the
Stasheff associahedron, has dimension $d-2$.)
These operations satisfy the curved (i.e., with $\mu^0$) $A_\infty$-relations
\vspace*{-2pt}\begin{equation}\label{eq:Ainfinity}
\sum_{\substack{0\leq k\leq d\\ 0\leq j\leq d-k}} (-1)^{*}
\mu^{d+1-k}(x_d,\dots,x_{j+k+1},\mu^k(x_{j+k},\dots,x_{j+1}),x_j,\dots,x_1)=0,
\vspace*{-3pt}\end{equation}
where $*=\deg(x_1)+\dots+\deg(x_j)-j$.  These relations express the fact
that the boundaries of 1-dimensional moduli spaces of solutions to Floer's
equation on the disc correspond to pairs of discs with the output of one matching an input of the
other; these can arise from nodal degenerations of the domain at the
boundary of the Stasheff associahedron (the terms in \eqref{eq:Ainfinity}
with $2\leq k\leq d-1$), breaking off of a Floer strip at an input ($k=1$) or
at the output ($k=d$), or disc bubbling ($k=0$).
When $\mu^0=0$, the first few relations (for $d=1,2,3$) state respectively that the differential
$\partial=\mu^1$ squares to zero, the product $\mu^2$ satisfies the Leibniz
rule, and $\mu^2$ is associative up to an explicit homotopy given by $\mu^3$.

The Floer differential generally fails to square to zero for Lagrangian
submanifolds with $\mu^0\neq 0$. Bounding
cochains were introduced by Fukaya et al.\ \cite{FO3book} as a way to
algebraically deform the $A_\infty$-operations to cancel $\mu^0$. 
Working over the Novikov field, an element $b\in CF^*(L,L)$ (of odd degree, and of positive Novikov
valuation, i.e.\ with coefficients only involving positive powers of the
formal variable $T$) is a {\bf bounding cochain} for $L$ if
\vspace*{-2pt}$$\mu^0_b = \sum_{k\geq 0} \mu^k(b^{\otimes k})=\mu^0 + \mu^1(b)+\mu^2(b,b)+\dots =
0\in CF^*(L,L).\vspace*{-4pt}$$
A Lagrangian submanifold $L$ is said to be {\bf unobstructed} if it admits
a bounding cochain (or one says that the pair $(L,b)$ is unobstructed).
More generally, a Lagrangian is {\em weakly
unobstructed} if it admits a
{\em weak bounding cochain}, i.e.\ $b\in CF^*(L,L)$ such that $\mu^0_b$
is a scalar multiple $\lambda\,1_L$ of the unit in $CF^*(L,L)$ for some
$\lambda\in\mathbb{K}$.

Given Lagrangian submanifolds $L_0,\dots,L_d$ and bounding cochains $b_i\in CF^*(L_i,L_i)$ for
$i=0,\dots,d$, one can modify the Floer operations \eqref{eq:mu_d} to
\vspace*{-4pt}\begin{equation}\label{eq:mu_d_b}\mu^d_b(x_d,\dots,x_1)=\sum_{k_0,\dots,k_d\geq 0\!\!}
\mu^{d+k_0+\dots+k_d}\bigl(b_d^{\otimes k_d},x_d,\dots,b_1^{\otimes
k_1},x_1,b_0^{\otimes k_0}\bigr).\vspace*{-4pt}\end{equation}
These modified operations satisfy the uncurved (i.e., without $\mu^0$)
$A_\infty$-relations; in particular $(\mu^1_b)^2=0$, so the Floer cohomology
$HF^*((L_0,b_0),(L_1,b_1))$ of a pair of Lagrangian submanifolds
equipped with bounding cochains is well-defined.
Likewise for weakly unobstructed Lagrangians that share a common value of the
constant $\lambda\in\mathbb{K}$.

\begin{remark}
There are several different possible ways to define the Floer complex $CF^*(L,L)$
of a compact Lagrangian submanifold with itself. Here we treat this case
identically to the Floer complex of a pair of distinct Lagrangians, by picking
a Hamiltonian perturbation (chosen generically so that $\varphi_H^1(L)$
meets $L$ transversely); this is similar to e.g.\ the construction in
\cite{SeBook}. However, there exist other models where $CF^*(L,L)$ is
defined to consist of singular chains or differential forms on $L$, and the
corresponding inputs or outputs of Floer operations are not boundary
punctures of the domain $S$ but rather boundary marked points,
with input and output data pulled back and pushed forward through evaluation
maps at these marked points; this is the approach taken in e.g.\
\cite{ChoOh}, \cite{FO3book}, \cite{SolomonTukachinsky}, etc.
Yet another option is to pick a Morse function $f:L\to\R$ and define
$CF^*(L,L)$ to be the free module generated by the critical points of $f$;
the Floer operations
then count (perturbed) ``treed $J$-holomorphic discs'', i.e.\ configurations of
$J$-holomorphic discs with boundary on $L$, together with
gradient flow lines of $f$ connecting their boundaries 
to each other and to the input and output critical points; see e.g.\
\cite{BiranCornea}, \cite{CharestWoodward}.
Each approach
has specific advantages and drawbacks; in the end they yield isomorphic
Floer cohomology algebras.
\end{remark}

\subsection{The Fukaya category.}

There are many versions of the Fukaya category in the literature, depending
on what kinds of Lagrangian submanifolds are allowed, whether a $\Z$-grading
is desired (this typically requires $2c_1(TM)=0$), the possible presence
of additional data such as bounding cochains and local systems, etc. In all
cases, the goal is to associate to a symplectic manifold an 
$A_\infty$-category $\mathcal{F}(M)$ whose objects are Lagrangian
submanifolds (together with additional data), with morphisms given by Floer complexes
and compositions given by the Floer operations \eqref{eq:mu_d}.
See e.g.\ \cite{AuBeginner} for an overview of the subject.

The simplest version of the Fukaya category, which can be defined over
any coefficient field, considers compact exact Lagrangian submanifolds in
a Liouville manifold $(M,\omega=d\lambda)$. In this case, exactness
precludes bubbling, so that unobstructedness holds tautologically without
the need to introduce bounding cochains, and a priori estimates on the
energy of Floer solutions (independently of homotopy class) make it
unnecessary to work over Novikov coefficients. (However, an object of
the Fukaya category should still come equipped with a spin structure,
so as to orient moduli spaces of Floer solutions, and grading data if
desired.) See \cite{SeBook} for a detailed treatment.

Outside of the exact setting, one should restrict oneself to unobstructed
(or weakly unobstructed) Lagrangian submanifolds, possibly after equipping
them with (weak) bounding cochains. It is also in general necessary to work over a Novikov
field and count solutions of Floer's equation with weights $T^{E(u)}$ 
determined by their energy \eqref{eq:energy}. 
Gromov compactness implies that the weighted counts of Floer solutions are
well-defined as elements of the Novikov field, even when
they involve contributions from 
infinitely many homotopy classes of discs.
(An exception is the case of {\em Bohr-Sommerfeld} or {\em balanced}
Lagrangian submanifolds in a monotone symplectic manifold, where the energy
of discs is controlled by their index.)

In some settings, e.g.\ for homological
mirror symmetry, one may also consider Lagrangian submanifolds equipped
with {\em local systems}. When working over the Novikov field, these
are required to be {\em unitary}, i.e.\ their holonomy
should be of the form $a_0+\sum a_i T^{\lambda_i}$ with
$a_0$ invertible and $\lambda_i>0$ for all $i$.
Solutions to
Floer's equation are then counted with weights
$w([u])=T^{E(u)}\,\mathrm{hol}(\partial u)$, where $E(u)$ is the energy 
and $\mathrm{hol}(\partial u)$ is the product of the holonomies
of the local systems along the boundary of the disc (see e.g.\ \cite[Remark
2.11]{AuBeginner}).
\medskip

Perhaps more importantly from a geometric perspective, there are a number of different versions of the Fukaya category 
for noncompact Lagrangian submanifolds in noncompact symplectic manifolds, 
depending on the conditions imposed on the geometric behavior of the
Lagrangians at infinity and on the class of Hamiltonian perturbations used
to define Floer complexes. 
The {\bf wrapped Fukaya category} of a Liouville
manifold involves exact Lagrangian submanifolds which are cylindrical at
infinity, with morphisms given by wrapped Floer complexes, i.e.\ the
direct limits of Floer complexes with respect to linear-growth Hamiltonian perturbations of 
increasingly large slope at infinity, and suitably constructed 
Floer $A_\infty$-operations (see e.g.\ \cite{AbGenerate, AS,GPS1} for
different approaches). There are also {\bf partially wrapped} Fukaya
categories \cite{Sylvan,GPS1,GPS2}, whose objects are required
to avoid certain directions at infinity (the ``stops''),
and morphisms are direct limits of Floer complexes with respect to
Hamiltonian perturbations whose flow stays away from the stops.
This includes Fukaya-Seidel
categories of symplectic fibrations \cite{SeBook,SeLF},
as well as other variants which lend themselves more easily to
calculations for certain classes of examples 
\cite{AuICM,Hanlon,AA} but can likely be translated into
the framework of \cite{GPS1,GPS2}.
One could also consider (partially) wrapped Floer theory on more general classes of
(not necessarily exact) noncompact symplectic manifolds. The simplest situation is
when the geometry at infinity reduces to the Liouville case (see e.g.\ \cite{RitterSmith}), but all that is required is some way to bound the
geometric behavior of Floer trajectories via maximum principles and/or energy estimates, 
and a suitable class of Hamiltonian perturbations; see e.g.\ \cite[Section~3]{AA}.
\medskip

The study of Kontsevich's homological mirror symmetry conjecture
\cite{KoICM} has led to a large body of work focused on computations of
Fukaya categories (ordinary, wrapped, or partially wrapped depending on the
geometric setting), with the goal of comparing these with categories of
coherent sheaves (or matrix factorizations) associated to mirror spaces.
While it is sometimes possible to list all objects of the Fukaya category
and directly compute their Floer cohomologies (see e.g.\ \cite{PZ}),
a more practical approach is usually to identify a collection of
objects which generate the Fukaya category. One says that a collection
of objects $\{G_i\}$ {\bf generates} (resp.\ {\bf split-generates}) an
$A_\infty$-category $\mathcal{C}$ if, in a triangulated enlargement of 
$\mathcal{C}$ (for instance twisted complexes or modules over $\mathcal{C}$
\cite{SeBook}), every object is quasi-isomorphic to an iterated mapping cone
built from (arbitrarily many copies of) the objects $G_i$ (resp.\ a direct summand in an iterated
mapping cone). Then $\mathcal{C}$ admits a fully faithful embedding into the category
of (right) modules over the $A_\infty$-algebra $\mathcal{A}=\bigoplus_{i,j}
\hom(G_i,G_j)$, the
{\em Yoneda embedding} taking each object $T$ to the $A_\infty$-module
$\bigoplus_j \hom( G_j,T)$ (equipped with structure maps
given by compositions in the category $\mathcal{C}$). In particular,
determining the Floer complexes of the (split-)generators and their
$A_\infty$ structure maps is sufficient to determine the whole Fukaya
category. See \cite{AuBeginner} for an informal treatment and \cite{SeBook}
for details. 

\begin{remark}
Floer-type homology invariants of 3-manifolds (resp.\ knots and links) often admit further 
categorifications, which associate an algebra to a surface (resp.\ a
configuration of points in the plane), and a module over this algebra 
to a 3-manifold with boundary (resp.\ a tangle). The existence of symplectic
interpretations of these invariants (cf.\ Remark \ref{rmk:lowdimtop})
suggests that the algebras of interest to low-dimensional topologists should
be understood as describing the Fukaya categories of the symplectic
manifolds that appear in this context. An example where this works well
is bordered Heegaard-Floer homology \cite{LOT}:
Lipshitz-Ozsv\'ath-Thurston's strands algebra is precisely the endomorphism
algebra of a collection of generators of the partially wrapped Fukaya
category of the symmetric product of a punctured surface, and the bordered
Floer modules associated to 3-manifolds with boundary can be viewed as
arising from the Yoneda embedding \cite{AuICM}. Another instance is the
arc algebra underlying Khovanov homology \cite{AbouzaidSmithKhovanov}.
\end{remark}

A number of generation criteria have been established to help
determine when certain objects generate the Fukaya category, starting with Seidel's
results on Fukaya categories in Lefschetz fibrations \cite{SeBook} and
Abouzaid's generation criterion for wrapped Fukaya categories
\cite{AbGenerate}, continuing with results on 
automatic generation in Calabi-Yau mirror symmetry
\cite{GPSh}, and generation results for (partially) wrapped Fukaya 
categories of Weinstein manifolds and sectors \cite{GPS2,CDRGG}.
Many of the examples we discuss below are in the context of homological
mirror symmetry, but first we start with an example that is of independent
interest to symplectic topologists.

\subsection{Cotangent bundles and the nearby Lagrangian conjecture.}
About twenty years ago, Nadler and Zaslow \cite{NZ,Nadler} constructed a quasi-equivalence between
the category of constructible sheaves
on a compact manifold $N$ and a certain ``unwrapped'' Fukaya
category of exact (not necessarily compact) Lagrangians in its cotangent
bundle $T^*N$, under which compact exact Lagrangian
submanifolds (with local systems) in $T^*N$ correspond to local systems on $N$.
Fukaya-Seidel-Smith independently arrived at essentially the same conclusion
using a different approach \cite{FSS}. This yields:

\begin{theorem}[Nadler, Fukaya-Seidel-Smith \cite{Nadler,FSS}]
\label{thm:nzfss}
Let $N$ be a compact spin manifold, and $L\subset T^*N$ be
a compact exact Lagrangian submanifold which is spin and whose 
Maslov class vanishes. Then $L$ is isomorphic to the zero section in
the Fukaya category of $T^*N$. In particular, the projection from $L$
to $N$ has degree $\pm 1$ and induces an isomorphism $H^*(L;\K)\simeq
H^*(N;\K)$ over any coefficient field.
\end{theorem}

Meanwhile, the structure of the wrapped Fukaya category of $T^*N$ was
elucidated by Abouzaid:

\begin{theorem}[Abouzaid \cite{AbCotangent}]
The wrapped Fukaya category of $T^*N$ is generated by a cotangent fiber
$F=T^*_qN$ (for any $q\in N$). In particular, $\W(T^*N)$ quasi-embeds into
the category of modules over $\mathrm{End}(F)\simeq C_{-*}(\Omega_q N)$,
the algebra of chains over the based loop space of $N$.
\end{theorem}

This is worth comparing with earlier work of
Viterbo, Abbondandolo-Schwarz and Salamon-Weber showing that the
symplectic cohomology $SH^*(T^*N)$ is isomorphic to the homology of the
free loop space $\mathcal{L}N$; see e.g.~\cite{Weber}.
(For simplicity we have stated all the above results under the assumption
that $N$ is spin; they continue to hold when $N$ is not spin if one
considers homology with suitably twisted coefficients.)

Returning to compact exact Lagrangians, we mention the following improvement
on Theorem \ref{thm:nzfss}:

\begin{theorem}[Abouzaid-Kragh \cite{AbKragh}]
Given any compact smooth manifold $N$, and any compact exact Lagrangian
submanifold $L\subset T^*N$, the projection from $L$ to $N$ is a
simple homotopy equivalence.
\end{theorem}

This is a significant partial result on
Arnold's {\bf nearby Lagrangian conjecture}, which asks
whether every compact exact Lagrangian submanifold in 
$T^*N$ is Hamiltonian isotopic to the zero section.
(By the Weinstein neighborhood theorem, a tubular neighborhood of a
Lagrangian submanifold $N\subset M$ is symplectomorphic to
a neighborhood of the zero section in $T^*N$, so Arnold's
conjecture indeed constrains nearby Lagrangians.)
Arnold's question remains open in general (though it has been answered positively in a
few cases), essentially because, even though Hamiltonian isotopic exact Lagrangian
submanifolds are Fukaya isomorphic, it is not clear that the converse should
hold. (Outside of the exact setting, there are infinite families of monotone 
Lagrangian tori in $\CP^2$ and many other symplectic manifolds which are
not Hamiltonian isotopic to each other yet define isomorphic objects of the
Fukaya category when equipped with suitable local systems; see e.g.\
\cite{Vianna}.)  In particular, Fukaya categories provide a wealth of
information about Lagrangian submanifolds in a given symplectic manifold, but they do not directly address
the problem of their classification up to Hamiltonian isotopy.
(That said, perhaps an even more basic question is whether objects
of the Fukaya category even behave in the geometric manner suggested by
the algebra; this is harder than one might think, already for surfaces
\cite{AurouxSmith}.)

That said, there is evidence that Floer theory can provide more refined
information, at least if one keeps track of more detailed information about
moduli spaces beyond mere counts of isolated solutions to Floer's equation
-- starting with the following result of Abouzaid:
\begin{theorem}[Abouzaid \cite{AbSpheres}]
Every homotopy sphere which embeds as a Lagrangian in $T^*S^{4k+1}$ bounds
a compact parallelizable manifold.
\end{theorem}
More recently, Floer homotopy theory has led to the introduction of ``spectral
Fukaya categories'', which are only beginning to be systematically explored
but already have applications to bordism and stable homotopy types of
quasi-isomorphic Lagrangians and to smooth structures on nearby Lagrangians;
see e.g.\ \cite{Large,PS1,PS2}.

\subsection{Open-closed maps, deformations, and homological mirror symmetry.}
\label{ss:ocdeforms}
Given Lagrangian submanifolds $L_0,\dots,L_d$ ($d\geq 0$) of
$(M,\omega)$,
moduli spaces of solutions to Floer's equation on discs with $d+1$ 
boundary punctures (as in Lagrangian Floer theory) and one single interior 
puncture (as in Hamiltonian Floer theory) can be used to define {\bf
open-closed} and {\bf closed-open} maps\vspace*{-1mm}
\begin{eqnarray}
\label{eq:OC}
&&OC^d:CF(L_{d},L_0)\otimes\dots\otimes CF^*(L_0,L_1)\to CF^*(M,H),\\
\label{eq:CO}
&&CO^d:CF^*(M,H)\otimes CF^*(L_{d-1},L_d)\otimes \dots\otimes
CF^*(L_0,L_1)\to CF^*(L_0,L_d),
\end{eqnarray}
\vspace*{-6mm}

\noindent
depending on whether the output of the operation corresponds to the 
interior puncture or a boundary puncture.
These maps intertwine the Floer differential on $CF^*(M,H)$ with the
differentials on the Hochschild complexes
\begin{eqnarray*}CC_*(\mathcal{F}(M),\mathcal{F}(M))&=&\bigoplus_{d\geq
0}\bigoplus_{L_0,\dots,L_d} CF^*(L_{d},L_0)\otimes\dots\otimes
CF^*(L_0,L_1) \qquad \text{and}\\
CC^*(\mathcal{F}(M),\mathcal{F}(M))&=&\prod_{d\geq
0}\prod_{L_0,\dots,L_d} \hom(CF^*(L_{d-1},L_d)\otimes\dots\otimes
CF^*(L_0,L_1),CF^*(L_0,L_d)),
\end{eqnarray*}
\vspace*{-4mm}

\noindent
thus inducing cohomology-level maps \vspace*{-2mm}
\begin{equation}\label{eq:OC_and_CO}
OC:HH_{*-\dim_\C M}(\F(M),\F(M))\to HF^*(M,H)
\quad \text{and}\quad CO:HF^*(M,H)\to HH^*(\F(M),\F(M)).
\end{equation}
In the compact setting, one can replace $HF^*(M,H)$ with the 
quantum cohomology $QH^*(M)$; and in the noncompact setting there
are similar maps in wrapped Floer theory, between the symplectic cohomology
$SH^*(M)$ and the Hochschild homology and cohomology of the wrapped Fukaya
category $\W(M)$. (However, open-closed maps work somewhat differently for partially wrapped Fukaya categories,
including Fukaya-Seidel categories.)

Open-closed maps play a key role in the structure of Fukaya categories. For
example:
\begin{theorem}[Abouzaid \cite{AbGenerate}, Ganatra \cite{Ganatra, GanatraCyclic}]
Let $(M,\omega)$ be a Liouville manifold, and assume there exists a
Hochschild cycle which maps to the unit $1\in SH^*(M)$ under the open-closed
map. Then:\vspace*{-2pt}
\begin{enumerate}
\item $\W=\W(M)$ is split-generated by the objects whose morphisms appear in the given
Hochschild cycle;\vspace*{-3pt}
\item $CO:SH^*(M)\to HH^*(\W,\W)$ is a ring isomorphism, and
$OC:HH_{*-\dim_\C M}(\W,\W)\to SH^*(M)$ is an isomorphism of modules over
$HH^*(\W,\W)\simeq SH^*(M)$;\vspace*{-3pt}
\item $\W$ is a homologically smooth (i.e., its diagonal bimodule is
perfect) non-compact Calabi-Yau category (i.e., the inverse dualizing
bimodule is isomorphic to the diagonal bimodule up to a grading shift);
\vspace*{-3pt}
\item the $S^1$-equivariant enhancement of $OC$ known as the cyclic open-closed map \cite{GanatraCyclic}
defines isomorphisms from the (positive, negative, periodic) cyclic homology
of $\W$ to (the corresponding variants of) the $S^1$-equivariant symplectic
cohomology of $M$.
\end{enumerate}
\end{theorem}
\begin{remark}
There is more to this story: the symplectic cohomology $SH^*(M)$ carries the further structure of a
homotopy BV-algebra, and the closed-open map is compatible with these
additional operations \cite{AGV,BEASh}.
\end{remark}
Similar results hold for Fukaya categories of compact symplectic
manifolds. These ingredients (and in particular the cyclic open-closed
map) play a key role in Ganatra-Perutz-Sheridan's program to recover
enumerative mirror symmetry from homological mirror symmetry for Calabi-Yau
varieties \cite{GPSh}. Namely, let $M,M^\vee$ be a pair of smooth projective
Calabi-Yau manifolds which satisfy homological
mirror symmetry in the sense that the Fukaya category $\F(M)$ is
quasi-equivalent to the category of coherent sheaves
of (a maximally unipotent degeneration of) $M^\vee$.
Passing to Hochschild cohomology, one obtains a ring isomorphism between
$HH^*(\F(M))\simeq QH^*(M)$ and $HH^*(\Coh(M^\vee))\simeq H^*(M^\vee,\bigwedge^*TM^\vee)$,
which is the starting point of ``classical'' mirror symmetry. Upgrading this
(as first proposed by Barannikov and Kontsevich) to an isomorphism of
variations of Hodge structures over formal punctured discs yields a
Hodge-theoretic version of mirror symmetry, which in turn has enumerative
consequences, such as e.g.\ the
classical mirror symmetry statement about counts of rational curves in the
quintic 3-fold \cite{GPSh}.
\medskip

Since Hochschild cohomology governs first-order deformations of the Fukaya
category, the closed-open map being an isomorphism implies that every
first-order deformation of the Fukaya category has a geometric origin.
In fact, a similar geometric construction also yields actual (formal) deformations
of the Fukaya category over the Novikov field. Given a class 
$\bb\in H^*(M,\K)$ (or more generally in Hamiltonian Floer cohomology) (of even degree, and of positive Novikov
valuation), the {\bf bulk-deformed} $A_\infty$ operations
$\mu^d_{\bb}=\mu^d+CO^d(\bb)+\dots$ are
defined by counting solutions to Floer's equation on a disc with $d+1$ boundary punctures
(or more if there are also bounding cochains) and any number $\ell\geq 0$
of (input) interior punctures, with (a chain-level representative of)
$\bb$ inserted at each of those. 
One then arrives at the bulk-deformed Fukaya category $\F(M,\omega,\bb)$. 
(Note that, since the curvature $\mu^0_\bb$ depends on $\bb$, so does the
notion of bounding cochain; for instance, for the purpose of proving
Theorem \ref{thm:floerlagFO3} it is enough to find a weak bounding cochain for
$L$ with respect to {\em some} bulk deformation.)

Another (related) situation where geometric deformations of
the Fukaya category naturally occur is when one compares the Fukaya category of a smooth
complex projective variety $M$ to that of the complement $M^0=M\setminus D$
of a smooth (or more generally, normal crossings) divisor $D\subset M$. Namely,
for Lagrangian submanifolds of $M$ which are disjoint from $D$, and using Floer
data that preserve positivity of intersections with $D$, one can
count solutions to Floer's trajectories in $M$, weighing those which
intersect the divisor $n$ times with an additional factor of $q^n$
(where $q$ is a formal variable). This gives rise to the {\bf relative
Fukaya category} $\F(M,D)$, which presents the Fukaya category of $M$
as a deformation of $\F(M^0)$ \cite{SeICM,PS}. The first order term
of this deformation, which corresponds to discs that intersect $D$
transversely once, can be expressed as the image under the closed-open map
of a certain element in $SH^*(M^0)$, sometimes called the
Borman-Sheridan class \cite{BEASh2}.

These considerations have led to an extremely effective approach to
homological mirror symmetry for vast classes of examples, starting with
Seidel's work on genus 2 surfaces \cite{SeGenus2} which start from 
a description as a
compactification of a $\Z/5$-cover of the pair of pants, continuing with
Sheridan's landmark result on Calabi-Yau hypersurfaces in projective space
\cite{Sheridan}, and more recently a proof of homological mirror symmetry
for Batyrev mirror pairs \cite{GHHPS}.
(That said, there are many other approaches to calculations of
Fukaya categories and proofs of homological
mirror symmetry; see e.g.\ \cite{LPpants}, \cite{AA} and \cite{Ueda} for some recent
examples.)

\begin{remark}
To explain the central role of relative Fukaya categories in many proofs
of homological mirror symmetry, recall that mirror
symmetry relates symplectic enumerative geometry near the {\em large volume limit}
(i.e., viewing formal variables
such as the Novikov parameter $T$ or the relative Fukaya category parameter
$q$ as living in a neighborhood of the origin) to algebraic geometry near the 
{\em large complex structure limit} (i.e.,
for varieties defined over formal power series, corresponding to maximally
degenerating families of complex varieties).
The complement of a sufficiently ample hypersurface (e.g.\ a
hyperplane section) in a projective variety is a Liouville (in fact, Stein)
manifold, whose Fukaya category can be determined by a wide range of
techniques, including microlocal sheaf theory or sectorial decompositions
(see the next section); the pushout diagrams inherent to sectorial descent
\cite{GPS2} for these large volume limits typically correspond under mirror
symmetry to gluing closed subschemes to produce a (singular, reducible)
large complex structure limit mirror (see e.g.\ \cite{GammageShende1,
GammageShende2}). One can then hope to recover mirror
symmetry for the projective variety itself by a deformation process as in
the above examples \cite{SeGenus2,Sheridan,GHHPS}, though the process for doing this in general has not been
fully elucidated.
\end{remark}

\begin{remark}
Deformation theory arguments also explain the appearance of {\bf
Landau-Ginzburg models}, i.e.\ spaces equipped with a function
called {\em superpotential}, and their categories of {\em matrix
factorizations}, in mirror symmetry outside of 
the Calabi-Yau setting. Namely, given a normal crossings anticanonical 
divisor $D$ in $M$ (i.e., Poincar\'e dual to $c_1(TM)$, which we assume to
be nef), 
the complement $M^0=M\setminus D$ is Calabi-Yau, and its (wrapped) Fukaya
category can often be related to coherent sheaves on a mirror space $M^\vee$.
Adding back in the divisor $D$ then deforms the Fukaya category as explained
above; however this deformation does not preserve the $\Z$-grading on
$\F(M^0)$, and the first-order deformation class 
lives in $HH^0(\F(M^0))$ rather than in degree
2 as in the Calabi-Yau case. Under mirror symmetry, this becomes an element
in $HH^0(\mathrm{Coh}(M^\vee))\simeq H^0(M^\vee,\mathcal{O}_{M^\vee})$,
i.e.\ a function on $M^\vee$: the superpotential.
Concretely, compact unobstructed Lagrangians in $M^0$ deform to
{\em weakly} unobstructed Lagrangians in $M$, with $\mu^0$ given by a
weighted count of holomorphic discs which intersect $D$ once: the (constant)
value of the superpotential on the support of the mirror object.
See also \cite{Au07}.
\end{remark}

\begin{remark}
Another important structural aspect of Fukaya-Floer theory that we haven't
touched on is its functoriality under Lagrangian correspondences: namely,
a Lagrangian submanifold in $(M\times N, -\pi_M^*\omega_M+\pi_N^*\omega_N)$,
or rather an object of its Fukaya category, determines an $A_\infty$-bimodule over the
Fukaya categories $\F(M)$ and $\F(N)$, often representable by an $A_\infty$-functor 
from $\F(M)$ to (an enlargement of) $\F(N)$. See e.g.\ \cite{MWW,Fcorr,AbBottman}.

There are also specific functors between the
Fukaya categories of a symplectic fibration (resp.\ sector),
its total space (resp.\ completion), and its fiber. Under mirror symmetry
these often correspond to inclusion and restriction of sheaves between a variety, a divisor in it, and
the divisor complement; see e.g.\ \cite[\S 5]{Au09} and
\cite{Sylvancapcup,NadlerLGCn,AuSpec,GammageJeffs,HanlonHicks}.
\end{remark}

\section{Towards the geometry of Floer theory.}

At this point, we have turned a lot of symplectic geometry into the language
of $A_\infty$ categories and homological algebra. This fits the philosophy
of ``noncommutative algebraic geometry'', i.e.\ the idea that any category 
that shares structural features with derived categories of algebraic varieties should be 
regarded as a noncommutative algebraic space. And yet, mirror symmetry tells
us that many Fukaya categories actually correspond to 
honest (commutative) algebraic spaces. This means that there should exist a (classical,
``commutative'') geometric perspective on Fukaya categories, 
lending itself to geometric decomposition results and local-to-global
principles -- the difficulty being that, a priori, Floer theory
is non-local, in the sense that even for Lagrangian submanifolds contained 
in a subset $U\subset M$ it may involve trajectories which do not remain within
$U$. In this section we discuss various perspectives that fit into this
general philosophy.

\subsection{Floer theory for families of Lagrangians.}\label{ss:family}

Consider a family of Lagrangian submanifolds $(F_b)_{b\in B}$ in a
symplectic manifold $(M,\omega)$, parametrized by a space $B$: for
instance a Lagrangian fibration $\pi:M\to B$
(possibly with singular fibers), but
one may also consider families in which
$F_b$ has different topology over different strata of $B$, etc.
Assuming the Lagrangians $F_b$ define objects of the Fukaya category of $M$,
one can associate to another object $L$ the collection of Floer 
complexes $CF^*(F_b,L)$, $b\in B$, and study how these depend on $b$.

For instance, Nadler and Zaslow's work on Fukaya categories of cotangent
bundles \cite{Nadler,NZ} fits into this philosophy. Consider the cotangent bundle $M=T^*N$ 
of a compact manifold
$N$, and the family of cotangent fibers $F_q=T^*_qN$, parametrized by $N$ itself.
Since nearby cotangent fibers are Hamiltonian isotopic to each other,
given an object $L$ of the Fukaya category supported on a compact exact
Lagrangian submanifold, a path from $q_0$ to $q_1$ in $N$ 
determines an isomorphism between $HF^*(F_{q_0},L)$ and $HF^*(F_{q_1},L)$,
so that the family of Floer cohomologies $(HF^*(F_q,L))_{q\in N}$ defines
a locally constant sheaf (a local system) over $N$. By the same process,
a non-compact exact
Lagrangian submanifold of $T^*N$, allowed to approach infinity along the 
conormal directions to some
given stratification of $N$, determines a constructible sheaf on $N$.

More generally, let $M$ be a {\bf Weinstein manifold} (i.e., a Liouville manifold
whose Liouville vector field is gradient-like; this includes Stein
manifolds, e.g.\ affine complex algebraic varieties). Then $M$ retracts
onto a half-dimensional {\em skeleton} (or {\em core}) $B$, whose strata are isotropic
submanifolds in $M$. When $B$ is sufficiently nice, one can view $M$ as
obtained by attaching the cotangent bundles of the critical
(top-dimensional, i.e.\ Lagrangian) strata of $B$ onto the subcritical 
strata (which do not
contribute to the Fukaya category) and onto each other. Taking Floer cohomology with
cocores (the analogues of cotangent fibers in this setting) turns Lagrangian
submanifolds of $M$ into sheaves on $B$ (to be suitably interpreted along the
singular locus of the skeleton) -- with the benefit that the new data is
local over $B$, in a way that wasn't a priori obvious for the Fukaya category 
(and especially for its wrapped variant). This ties in with a conjecture of Kontsevich
\cite{KoSheaf}, who proposed that the wrapped Fukaya category of the Weinstein manifold
$M$ should be given by global sections of a certain cosheaf of categories
over its skeleton. Nadler further proposed that the wrapped 
Fukaya category should be modelled by {\bf wrapped microlocal
sheaves} over the skeleton \cite{NadlerWrapped}. Due to the entirely
combinatorial / topological nature of microlocal sheaf theory, the latter
categories are (for trained experts at least) eminently more computable than
(wrapped or partially wrapped) Fukaya categories -- even as the
geometric aspects of wrapped Floer theory, and in particular the Hamiltonian dynamics of wrapping,
recede out of sight, since wrapped microlocal sheaves are defined in a purely
categorical manner as compact objects in the category of microlocal sheaves
\cite{NadlerWrapped}.  

Nadler's proposal has been proven by Ganatra,
Pardon and Shende, under the topological
assumption that $M$ admits a {\em stable polarization} (i.e., the direct sum
of $TM$ with a trivial vector bundle is the complexification of a real
vector bundle): namely, the wrapped Fukaya category of a stably polarized
Weinstein sector is equivalent to the category of wrapped microlocal sheaves
over its skeleton \cite{GPS3}. This validates a posteriori the approach to
homological mirror symmetry via calculations of categories of constructible sheaves pursued
by various authors during the 2010s (see e.g.\ \cite{FLTZ,NadlerLGCn}).
This is also the aproach used by Gammage and Shende to verify homological
mirror symmetry for hypersurfaces in $(\C^*)^n$, in the direction that
relates the wrapped Fukaya category of the hypersurface to coherent 
sheaves on a mirror space \cite{GammageShende1}. (The other direction,
relating coherent sheaves on the hypersurface to the Fukaya category of a
mirror, is proved in \cite{AA} using very different methods.)
\medskip

While the above discussion stems from considering
families of contractible Lagrangian submanifolds such as cotangent fibers, another setting in which
Floer theory for families has led to powerful advances is that
of {\bf Lagrangian torus fibrations}, which are at the heart of the Strominger-Yau-Zaslow
(SYZ) approach to mirror symmetry \cite{SYZ} (see also \cite{Au07}).
(Lagrangian torus fibrations arise in a variety of contexts, ranging from
completely integrable systems and other situations involving symplectic reduction,
to degenerations of algebraic varieties to unions of toric varieties.)
Given a Lagrangian torus fibration $\pi:M\to B$ (possibly with singular
fibers), rather than just the fibers $F_b=\pi^{-1}(b)$, one considers families
of objects of $\F(M)$ given by the Lagrangians $F_b$ together 
with choices of local systems and/or bounding cochains.
(One typically first considers the objects supported on the smooth fibers of
$\pi$, before possibly completing the family with objects
supported on or near the singular fibers.) When the fibers
$F_b$ are tautologically unobstructed (do not bound any 
holomorphic discs in $M$), the appropriate parameter space is the {\em
uncorrected}\/ {\bf SYZ mirror} $Y^0=\{(F_b,\xi)\}$, the space of pairs consisting of a smooth
fiber $F_b=\pi^{-1}(b)$ ($b\in B^0=B\setminus \mathrm{critval}(\pi)$),
together with a unitary rank 1 local system $\xi$ over $F_b$, which can be characterized by
its holonomy $\hol_\xi\in \hom(\pi_1(F_b), U(1)_\K)$. This carries
a ``dual'' torus fibration $\pi^\vee:Y^0\to B^0$ mapping $(F_b,\xi)$ to
$b$, and has a natural structure of analytic space over $\K$, for which the
Floer weights $z_\beta=T^{\smash{\int_\beta}\omega}\,
\hol_\xi(\partial\beta)$ of homotopy classes of discs $\beta\in
\pi_2(M,F_b)$ are locally analytic functions \cite{Au07,AbICM,Ffamily}.

In the presence of singular fibers and/or holomorphic discs, one should
``correct'' $Y^0$ to a moduli space $Y$ of weakly unobstructed objects of $\F(M)$
supported on the fibers of $\pi$ (the {\em corrected} SYZ mirror).
The fibers of $\pi$ which bound Maslov index 0 holomorphic discs typically concentrate
along a union of codimension 1 {\em walls} in $B$, across which the
statement that Hamiltonian isotopic Lagrangians define isomorphic objects
of the Fukaya category needs to be corrected by a modification of the
bounding cochain or of the local system \cite{FO3book}.  Regluing the pieces
of $Y^0$ over the chambers of $B$ via these wall-crossing transformations
yields the corrected SYZ mirror $Y$.  Outside of the Calabi-Yau setting,
the fibers of $\pi$ typically also bound Maslov
index 2 discs; these do not require a modification of the mirror geometry,
but they make the objects $(F_b,\xi)$ {\em weakly unobstructed}; the
quantity $W(F_b,\xi)$ such that $\mu^0=W(F_b,\xi)\,1_{(F_b,\xi)}$ defines
a global analytic function on $Y$, the {\em mirror superpotential}, thus making the
SYZ mirror a Landau-Ginzburg model.
See e.g.\ \cite{Au07,AAK,Tu,Yuan}.\smallskip

The {\bf family Floer program} initiated by Fukaya \cite{Ffamily} and
further developed by Abouzaid \cite{AbICM,AbFamily}, Tu \cite{Tu} and Yuan
\cite{Yuan} leverages Fukaya's
observation that the Floer complexes $CF^*((F_b,\xi),L)$ and their
differentials have a locally analytic dependence on $(F_b,\xi)\in Y$
\cite{Ffamily,Fadic} to construct the {\em family Floer functor}, from
the Fukaya category of $M$ to complexes of (coherent) analytic sheaves (or matrix
factorizations when the objects are only weakly unobstructed) on $Y$,
under the assumption that the fibration $\pi$ admits a Lagrangian
section. (Otherwise one obtains {\em twisted} sheaves with respect to some
gerbe on $Y$ \cite{AbICM}.)
Roughly speaking, the functor associates to an object $L$ of $\F(M)$ the
family of Floer complexes $CF^*((F_b,\xi),L)$, $(F_b,\xi)\in Y$, and to
a morphism $x\in CF^*(L_0,L_1)$ the family of maps from $CF^*((F_b,\xi),L_0)$
to $CF^*((F_b,\xi),L_1)$ given by Floer product with $x$. (In the weakly
unobstructed case, when $\mu^0_{(F_b,\xi)}=W(F_b,\xi)\,1$ and
$\mu^0_L=\lambda\,1$, the differential on the Floer complexes squares to
$(W-\lambda)\,\mathrm{id}$, hence the appearance of matrix factorizations.)
The family Floer functor was shown by Abouzaid to be a quasi-equivalence in
the case where $\pi:M\to B$ doesn't have any singular fibers, thereby
proving homological mirror symmetry for such SYZ mirror pairs \cite{AbFamily}.
To go beyond this, one needs to have a good understanding of homological mirror 
symmetry for neighborhoods of singular Lagrangians; this has now been done
explicitly for the local models that are relevant to SYZ mirror symmetry 
\cite{AbSylvan,GammageSYZ}. One should also note the ``localized mirror
functor'' approach of Cho-Hong-Lau \cite{CHL,CHL2}, which can in principle be
applied to very general singular Lagrangians and their nearby smoothings,
at the expense of involving noncommutative language.
\medskip

Another perspective on SYZ mirror symmetry (and conjecturally on family
Floer theory) is that proposed in \cite{negmaslov} to deal with the
noncommutative corrections to the mirror geometry that can arise for
pairs $(M,D)$ where the anticanonical divisor $D$ contains rational curves
of negative Chern number. The Floer theory of a smooth fiber
$F_b$ of the Lagrangian fibration $\pi:M\to B$, with coefficients
in a completion of the group ring $\K[H_1(F_b)]$ (counting discs with
universal holonomy weights), encodes the family of Floer cohomologies for
all the objects of $\F(M)$ supported on $F_b$ -- an idea already used to great
effect by Abouzaid in \cite{AbFamily}. The completions of
$\K[\pi_1(H_b)]$ assemble into a sheaf $\O_{an}$ over the smooth locus $B^0$, 
namely the pushforward of the structure sheaf of the uncorrected SYZ mirror
$Y^0$ under the projection $\pi^\vee:Y^0\to B^0$.  Moduli spaces of 
pseudo-holomorphic discs in $M$ with boundary in the
fibers $F_b$ (where $b$ is now allowed to vary over $B^0$) then determine 
$A_\infty$-operations $\{\mu^k\}_{k\geq 0}$ not just on Floer cochains of a fixed fiber $F_b$
with coefficients in $\K[H_1(F_b)]$, but also on cochains on
$M^0=\pi^{-1}(B^0)\subset M$ with
coefficients in the pullback of $\O_{an}$, or equivalently, on
$\mathfrak{C}=C^*(B^0\,;\,C^*(F_b)\otimes\O_{an})$. The curvature $\mu^0$ of
this $A_\infty$-algebra encodes in principle all the information needed to
understand corrections to the mirror geometry \cite{negmaslov}. 

A conjectural systematic approach to these is as follows.
$H^*(F_b)\otimes \K[H_1(F_b)]$ is isomorphic to the homology of the free loop space
$\mathcal{L}F_b$, and carries a natural Lie bracket of degree
$-1$ (the Chas-Sullivan string bracket),\vspace*{-4pt}
\begin{equation}\label{eq:HF-bracket}
\{z^\gamma\, \alpha, z^{\gamma'} \alpha'\}=z^{\gamma+\gamma'}\,
\bigl(\alpha\wedge (\iota_\gamma \alpha')+(-1)^{|\alpha|}
(\iota_{\gamma'}\alpha)\wedge \alpha'\bigr).\vspace*{-4pt}
\end{equation}
Since $H^*(F_b,\R)\simeq \Lambda^* H^1(F_b,\R)\simeq\Lambda^* TB$, 
there is a natural map from $H^*(F_b)\otimes \K[H_1(F_b)]$ to polyvector
fields over (a neighborhood of) $(\pi^\vee)^{-1}(b)\subset Y^0$, under which
$\{\cdot,\cdot\}$ corresponds to the Schouten-Nijenhuis bracket.
It is conjectured in \cite{negmaslov} that, for a suitable model of Floer
theory, the curvature $\mu^0$ of $\mathfrak{C}$ satisfies the Maurer-Cartan
equation $\delta\mu^0+\frac12\{\mu^0,\mu^0\}=0$ (or its $L_\infty$ analogue)
with respect to the classical differential and the bracket \eqref{eq:HF-bracket}. 
The corrected differential
$\delta+\{\mu^0,\cdot\}$ on $\mathfrak{C}$ then squares to zero (unlike
$\mu^1$). Conjecturally, under mirror symmetry this
amounts to deforming
the \v{C}ech complex of polyvector fields $C^*(Y^0;\Lambda^*
TY^0)$ to arrive at polyvector fields (or their appropriate noncommutative
analogue) on the corrected mirror.
Moreover, this differential is expected to be part of a homotopy
BV-algebra structure on $\mathfrak{C}$, which should match that on polyvector
fields (or Hochschild cochains) under mirror symmetry.
This story should also be compatible with the construction of family Floer
functors; a first step can be found in Hoek's thesis \cite{Hoek}, which constructs a
functor from Lagrangian sections of $\pi$ to modules over the
curved $A_\infty$-algebra $(\mathfrak{C},\{\mu^k\})$.

\subsection{Geometric decompositions and local-to-global principles.}

The two constructions discussed in the previous section turn objects of the Fukaya 
category of a symplectic manifold $M$ into sheaves on another space (the skeleton 
$B$ in the Weinstein setting, the SYZ mirror space $Y$ in the setting of Lagrangian
torus fibrations), which have locality properties. Despite the non-local
nature of Lagrangian Floer theory, in both cases one rightfully expects
geometric decompositions of $M$ into suitable pieces to correspond to
geometric decompositions of $B$ and $Y$, thus suggesting that the Fukaya
category of $M$ can be computed from those of the pieces.  There are two
very different manners in which this principle works in practice.

\subsubsection*{Covariant functoriality: Liouville sectors.}

Liouville sectors, introduced and studied by Ganatra, Pardon and Shende
\cite{GPS1,GPS2}, are a class of exact symplectic manifolds with boundary, whose
non-compact ends are convex, i.e.\ modelled on the positive symplectization 
of a contact manifold with boundary, and which satisfy a set of geometric conditions
along their boundary ensuring that (for suitable almost-complex
structures) families of $J$-holomorphic curves cannot escape through the
sectorial boundary. A simple example of sector is the cotangent bundle
$T^*N$ of a manifold with boundary $N$; in this case the sectorial boundary is $\partial(T^*N)\simeq \R\times
T^*(\partial N)$. More generally, {\em Weinstein sectors} are built by
attaching cotangent bundles of critical strata of their skeleton 
(along part of their boundary only) onto a subscritical sector; and a
Liouville manifold can be decomposed along {\em sectorial hypersurfaces}
into a union of sectorial pieces.

The objects of the (partially) wrapped
Fukaya category of a Liouville sector \cite{GPS1} are properly embedded
exact Lagrangian submanifolds which are cylindrical at infinity and stay away from the sectorial
boundary;
morphisms are given by homotopy colimits of Floer complexes for the
images of the Lagrangians under
increasingly large Hamiltonian perturbations in the Reeb direction
at infinity (without crossing the sectorial
boundary). The key result about the behavior of these categories is:

\begin{theorem}[Ganatra-Pardon-Shende \cite{GPS1,GPS2}]\label{thm:gps}
Wrapped Fukaya categories of Liouville sectors are covar\-iantly functorial,
i.e.\ a proper inclusion $i:M'\hookrightarrow M$ gives rise to a functor
$i_*:\W(M')\to W(M)$. Moreover, when a Liouville manifold (or
sector) $M=M_1\cup_\partial M_2$ is obtained by gluing two sectors $M_1,M_2$ along
the common sectorial
boundary $\partial M_i\simeq \R\times F$, the inclusion functors give rise
to a pushout diagram
\begin{equation}\label{eq:gps}
\begin{CD}
\W(F)@>>> \W(M_1)\\
@VVV @VVV\\
\W(M_2)@>>> \W(M)
\end{CD}
\end{equation}
\end{theorem}
(Note that $\W(F)$ is equivalent to the wrapped
Fukaya category of the {\em stabilized} sector $T^*[0,1]\times F$, which admits
proper inclusions to neighborhoods of the sectorial boundaries of $M_1$ and
$M_2$.) This allows one to compute $\W(M)$ as the colimit of a diagram
involving the inclusion functors from $\W(F)$ into $\W(M_1)$ and $\W(M_2)$.
(Similarly for decompositions of a Liouville manifold or sector into more than two sectorial pieces.)
See also \cite{GanatraICM} for further discussion of Liouville sectors and
their Fukaya categories.

The gluing of wrapped Fukaya categories along sectorial decompositions
typically corresponds under homological mirror symmetry to the gluing of
categories of coherent sheaves under decompositions into unions of closed
subschemes. 

\begin{example} A twice-punctured torus $M$ can be cut along two arcs into a union of
two sectors $M_1,M_2$, both isomorphic to a cylinder with two stops (boundary
arcs). The corresponding inclusion functors give rise to a pushout diagram of 
wrapped Fukaya categories which, under homological mirror symmetry, matches 
the diagram of inclusion functors for coherent sheaves on the union of two
copies of the projective line, $Y_1,Y_2\simeq \PP^1$, glued together at two points (say
$\{0,\infty\}\subset Y_i$), to produce a twice nodal elliptic curve
$Y=Y_1\cup_{\{0,\infty\}}Y_2$.
\medskip

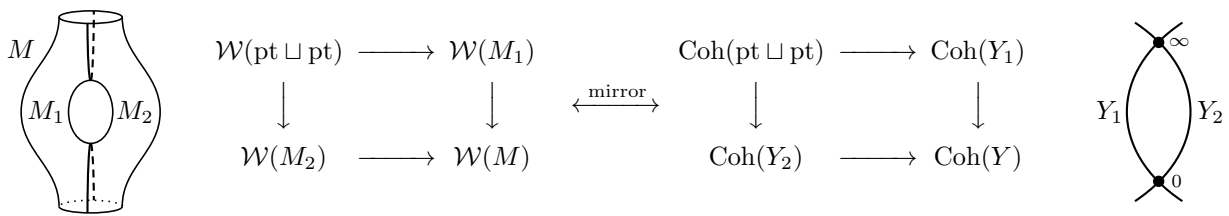
\centerline{
\begin{tikzpicture}[scale = 0.42, baseline=(current bounding box.center)]
\draw[semithick] (0,3) ellipse (1 and 0.2);
\draw[semithick](-1,3) to [out=-90,in=90](-2.2,0) to [out=-90,in=90] (-1,-3);
\draw[semithick](1,3) to [out=-90,in=90](2.2,0) to [out=-90,in=90] (1,-3);
\draw[semithick,dotted] (1,-3) arc (0:180:1 and 0.2);
\draw[semithick] (-1,-3) arc (180:360:1 and 0.2);
\draw[semithick](0,0) ellipse (0.7 and 1);
\draw[thick](-0.1,2.8) arc (180:270:0.1 and 1.8);
\draw[thick,densely dashed](0,1) arc (-90:0:0.1 and 2.2);
\draw[thick](-0.1,-3.2) arc (180:90:0.1 and 2.2);
\draw[thick,densely dashed](0,-1) arc (90:0:0.1 and 1.8);
\node at (-1.4,0){$M_1$};
\node at (1.4,0){$M_2$};
\node at (-2.2,2){$M$};
\end{tikzpicture}
\qquad
$\begin{CD}
\!\!\W(\mathrm{pt}\sqcup\mathrm{pt})@>>> \W(M_1)\\
@VVV @VVV\\
\W(M_2)@>>> \W(M)
\end{CD}$
\quad $\xleftrightarrow{\text{\ mirror\ }}$\quad
$\begin{CD}
\!\!\Coh(\mathrm{pt}\sqcup\mathrm{pt})@>>> \Coh(Y_1)\\
@VVV @VVV\\
\Coh(Y_2)@>>> \Coh(Y)
\end{CD}$
\qquad
\begin{tikzpicture}[scale = 0.42, baseline=(current bounding box.center)]
\draw[thick](0.75,2.828) arc (135:225:6 and 4);
\draw[thick](-0.75,2.828) arc (45:-45:6 and 4);
\fill(0,2.2) circle (0.18);
\fill(0,-2.2) circle (0.18);
\node at (0.65,2.2){$\scriptstyle \infty$};
\node at (0.57,-2.2){$\scriptstyle 0$};
\node at (-1.55,0){$Y_1$};
\node at (1.65,0){$Y_2$};
\end{tikzpicture}
}
\end{example}

This perspective is particularly useful for studying mirror symmetry at
the large limit, e.g.\ it can be used to understand Gammage and Shende's
result on mirror symmetry for very affine hypersurfaces \cite{GammageShende1}.

\subsubsection*{Contravariant functoriality.}

The above covariant functoriality notwithstanding, wrapped Fukaya categories
generally behave in a contravariant
manner with respect to inclusions of subdomains. 
The prototypical instance of this is 
the {\em restriction functor} constructed
by Abouzaid and Seidel \cite{AS} from the wrapped Fukaya category of a
Liouville manifold $M$ to that of a {\em Liouville subdomain} $M_{-}\subset M$
(i.e., a codimension 0 submanifold with convex boundary transverse to the 
Liouville vector field; the wrapped Fukaya category of $M_-$ is defined to
be that of its {\em completion}, the Liouville manifold
$\hat{M}_-=M_{-}\cup_{\partial M_-}\,[1,\infty)\!\times\! \partial M_{-}$).
The construction relies on the same geometric idea as Viterbo's earlier
construction of a restriction map 
from the symplectic cohomology of $M$ to that of $M_-$ \cite{Viterbo},
namely the use of Hamiltonians which grow very steeply near $\partial
M_-$, for which there are no Floer trajectories from generators in $M\setminus M_-$ to 
generators in $M_-$, so that the chain-level restriction map amounts to a
quotient by a subcomplex (in fact, an ideal with respect to the product operations).

Heather Lee's thesis \cite{Lee} used a similar idea to construct
restriction functors associated to decompositions of Riemann surfaces into
pairs of pants (overlapping in cylinders) and show that these give rise to
pullback diagrams. (While Lee worked in the exact setting, the result also
applies to
Fukaya categories of closed surfaces \cite{AurouxSmith}.)

\begin{theorem}[Lee \cite{Lee}]\label{thm:lee}
Given a decomposition of a Riemann surface $\Sigma$ into a union of two
subsurfaces $\Sigma_1\cup_\gamma \Sigma_2$ obtained by cutting $\Sigma$ along a 
nontrivial simple closed curve (or union of such curves) $\gamma$, there are restriction functors from the
wrapped Fukaya category of $\Sigma$ to those of $\Sigma_1,\Sigma_2$ and 
to a cylindrical neighborhood $C$ of $\gamma$,
which form a pullback diagram
\begin{equation}\label{eq:hlee}
\begin{CD}
\W(\Sigma)@>>> \W(\Sigma_1)\\
@VVV @VVV\\
\W(\Sigma_2)@>>> \W(C)
\end{CD}
\end{equation}
\end{theorem}

\noindent (As before, by $\W(\Sigma_i)$ and $\W(C)$ we mean the wrapped Fukaya categories of the
completions of these subsurfaces.)

Lee further showed that, for pair of pants decompositions of complex curves
in $(\C^*)^2$ near the tropical limit, the diagram \eqref{eq:hlee}
corresponds under mirror symmetry to the restriction functors associated to
a natural affine cover of the mirror, i.e.\ the computation of wrapped Floer
cohomology for a pair of Lagrangians in $\Sigma$ via the pullback diagram
matches the computation of sheaf cohomology via \v{C}ech complexes on
the mirror \cite{Lee}.

See also \cite{PascaleffSibilla} for a different perspective which leads to
very similar results, but allows for further generalizations and
applications to homological mirror symmetry in higher dimensions \cite{MSZ}.

\begin{example} As above, we consider a twice-punctured torus, 
but now we decompose it along two simple closed curves into two pairs of
pants
$\Sigma_1,\Sigma_2$, glued to each other along two
cylindrical necks $C',C''$. The pair of pants is mirror to the affine curve
$\{xy=0\}\subset \mathbb{A}^2$ (cf.\ e.g.\ \cite[\S 3]{AuSpec}), while
the cylinder is mirror to the punctured affine line $\mathbb{G}_m=\mathbb{A}^1-\{0\}$,
and the restriction functors on wrapped Fukaya categories correspond under
mirror symmetry to the restriction functors associated to the decomposition of $Y=\PP^1\cup_{\{0,\infty\}} \PP^1$
into two affine charts $U_1,U_2$ (the complement of $0$ and the complement of
$\infty$), with overlap $U_1\cap U_2\simeq \mathbb{G}_m\sqcup \mathbb{G}_m$.
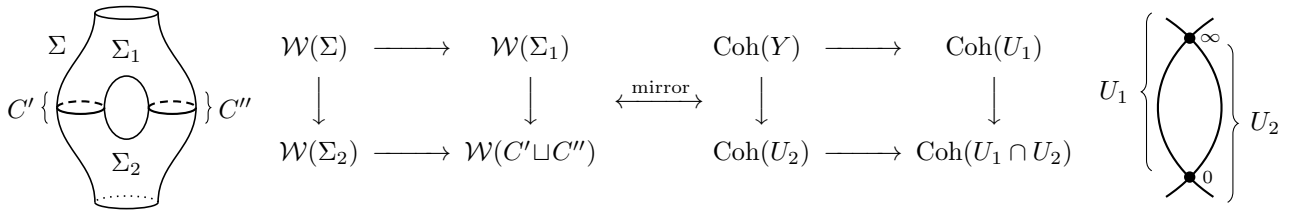
\medskip

\centerline{
\begin{tikzpicture}[scale = 0.42, baseline=(current bounding box.center)]
\draw[semithick] (0,3) ellipse (1 and 0.2);
\draw[semithick](-1,3) to [out=-90,in=90](-2.2,0) to [out=-90,in=90] (-1,-3);
\draw[semithick](1,3) to [out=-90,in=90](2.2,0) to [out=-90,in=90] (1,-3);
\draw[semithick,dotted] (1,-3) arc (0:180:1 and 0.2);
\draw[semithick] (-1,-3) arc (180:360:1 and 0.2);
\draw[semithick](0,0) ellipse (0.7 and 1);
\draw[thick](-2.2,0) arc (180:360:0.75 and 0.2);
\draw[thick](0.7,0) arc (180:360:0.75 and 0.2);
\draw[thick, densely dashed](-2.2,0) arc (180:0:0.75 and 0.2);
\draw[thick, densely dashed](0.7,0) arc (180:0:0.75 and 0.2);
\draw[decorate, decoration={brace, amplitude=2pt}](-2.5,-0.5)--(-2.5,0.5) 
 node[midway,xshift=-0.35cm]{$C'$};
\draw[decorate, decoration={brace, amplitude=2pt}](2.5,0.5)--(2.5,-0.5) 
 node[midway,xshift=0.4cm]{$C''$};
\node at (0,1.8){$\Sigma_1$};
\node at (0,-1.8){$\Sigma_2$};
\node at (-2.2,2){$\Sigma$};
\end{tikzpicture}
\ %
$\begin{CD}
\!\!\W(\Sigma)@>>> \W(\Sigma_1)\\
@VVV @VVV\\
\W(\Sigma_2)@>>> \W(C'\!\sqcup\! C'')
\end{CD}$
\ $\xleftrightarrow{\text{\ mirror\ }}$\ %
$\begin{CD}
\!\!\Coh(Y)@>>> \Coh(U_1)\\
@VVV @VVV\\
\Coh(U_2)@>>> \Coh(U_1\cap U_2)
\end{CD}$
\ %
\begin{tikzpicture}[scale = 0.42, baseline=(current bounding box.center)]
\draw[thick](0.75,2.828) arc (135:225:6 and 4);
\draw[thick](-0.75,2.828) arc (45:-45:6 and 4);
\fill(0,2.2) circle (0.18);
\fill(0,-2.2) circle (0.18);
\node at (0.65,2.2){$\scriptstyle \infty$};
\node at (0.57,-2.2){$\scriptstyle 0$};
\draw[decorate, decoration={brace, amplitude=4pt}](-1.2,-2)--(-1.2,3) 
 node[midway,xshift=-0.5cm]{$U_1$};
\draw[decorate, decoration={brace, amplitude=4pt}](1.2,2)--(1.2,-3) 
 node[midway,xshift=0.5cm]{$U_2$};
\end{tikzpicture}
}
\medskip
\end{example}

It is important to note that Theorem \ref{thm:lee} fails when one of the
curves is homotopically trivial, i.e.\ one of the pieces of the
decomposition is a disc. Indeed, in this case the locality of Floer
trajectories for sufficiently steep Hamiltonians fails to hold:
capping $\Sigma_i$ with a disc deforms its Fukaya category nontrivially
(cf.\ \S \ref{ss:ocdeforms}), and any
local-to-global principle in such a setting must incorporate this
deformation.

In the ``closed-string'' setting, an important conceptual step in this direction
is the work of Groman and Varolgunes \cite{GV} defining the {\em relative
symplectic cohomology} $SH^*_M(V)$ of a domain $V\subset M$, which 
deforms $SH^*(V)$ by the contributions of portions of Floer solutions in
$M\setminus V$ connecting generators in $V$. (When $V$ is the complement of
a divisor $D$ in $M$, this amounts to the description of Floer theory
in $M$ as a deformation of Floer theory in $M\setminus D$, but the setting
is much more general, allowing one to consider e.g.\ subsets of the form
$\pi^{-1}(U)$ for suitable $U\subset B$ in the total space of a Lagrangian fibration $\pi:M\to B$.)
One expects that a similar notion should exist for wrapped Fukaya categories 
and give rise to pullback diagrams of restriction functors associated
to covers by subdomains, even in settings where the local categories have to be deformed by
non-local contributions.

In the setting of SYZ mirror symmetry for a Lagrangian torus fibration
$\pi:M\to B$, the relative symplectic cohomologies and relative wrapped Fukaya
categories should behave like sheaves with respect to 
a cover of $B$ by suitable domains. Returning to the discussion at 
the end of \S \ref{ss:family}, the same should also be true of the
corrected family Floer complex $(\mathfrak{C},\delta+\{\mu^0,\cdot\})$ 
introduced in \cite{negmaslov}; in fact the latter should be
quasi-isomorphic (compatibly with $BV$-algebra structures) to the \v{C}ech complex
of relative symplectic cohomologies, since both are expected to give \v{C}ech models for
polyvector fields or Hochschild cochains on the corrected mirror.


\end{document}